\documentclass{conm-p-l}
\usepackage{epsfig}
\usepackage{amsmath,amssymb,amsthm}

\begin{document}

\theoremstyle{definition}

\newtheorem{Lemma}{Lemma}[section]
\newtheorem{Theorem}[Lemma]{Theorem}
\newtheorem{Proposition}[Lemma]{Proposition}
\newtheorem{Corollary}[Lemma]{Corollary}
\newtheorem{Remark}[Lemma]{Remark}
\newtheorem{Definition}[Lemma]{Definition}
\newtheorem{Question}[Lemma]{Question}
\newtheorem{theorem}{Theorem}[section]
\newtheorem{lemma}[theorem]{Lemma}

\numberwithin{equation}{section}

\newcommand{\beqn}{\begin{equation}}
\newcommand{\eeqn}{\end{equation}}
\newcommand{\nn}{\nonumber}
\newcommand{\la}{\langle}
\newcommand{\ra}{\rangle}
\newcommand{\cleq}{{\preccurlyeq}}
\def\R{{\mathbb R}}
\def\C{{\mathbb{C}}}
\def\pa {{\partial}}
\def\ep{{\epsilon}}
\newcommand{\ve}{\varepsilon}

\def\bn{{\bf n}}
\newcommand{\grad}{{\nabla} }

\title{The point source inverse back-scattering problem}

\author{Rakesh}
\address{Department of Mathematical Sciences, University of Delaware, Newark, DE 19716, USA}
\email{rakesh@math.udel.edu}
\thanks{The first author was partially supported by NSF grants DMS 0907909, DMS 1312708.}

\author{Gunther Uhlmann}
\address{Department of Mathematics, University of Washington, Seattle, WA 98195, USA}
\email{gunther@math.washington.edu}
\thanks{The second author was partially supported by the NSF and a Simons Fellowship.}

\subjclass[2010]{35R30, 35P25}
\date{February 24, 2014}

\dedicatory{Dedicated to D.H. Phong on the occasion of his 60th birthday}

\keywords{inverse backscattering, inverse point source, spherical transform}

\begin{abstract}
We consider the inverse problem of recovering a potential by measuring the response at a point to a source located at the same point and then varying the point on the surface of a sphere. This is a similar to the inverse backscattering problem. We show that if the angular derivatives of the difference of two potentials having the same data is controlled by the $L^2$ norm of the difference of the potentials they must be equal. In particular this shows injectivity of the inverse problem for radial potentials.
\end{abstract}

\maketitle


\section{Introduction}

A long standing formally determined inverse problem is the inverse back-scattering problem - see [RU14] and the references there for 
details. In this article we study a similar formally determined problem except we use point sources instead of plane waves and the 
given data is, the response of the medium measured at the same location as the point source, instead of the far field pattern as in the 
standard inverse back-scattering problem. Below $B$ denotes the {\bf open} unit ball in $R^3$ and $S$ the unit sphere in $\R^3$.

Let $q(x)$ be a smooth function on $\R^3$ which is supported in $B$. For each $a \in S$, let $U^a(x,t)$ be the solution of
 the IVP
\begin{align}
(\Box - q)U^a(x,t) = \delta(x-a,t) & \qquad 
(x,t) \in \R^3 \times \R 
\label{eq:Ude}\\
U^a(x,t)=0 & \qquad \mbox{for} ~ t<0 ~.
\label{eq:Uic}
\end{align}
The goal is the recovery of $q(.)$ if $U^a(a,t)$ is known for
all $a \in S$ and all $t \in (0,2]$.

As shown in [Rom74], we may express $U^a(x,t)$ as
\[
U^a(x,t) = \frac{1}{4 \pi} \frac{ \delta( t- |x-a|)}{|x-a|} + u^a(x,t)
\]
where $u^a(x,t)$ is a function which is zero on $t< |x-a|$ and, in the region $t \geq |x-a|$, $u^a$ is the smooth solution of the 
Goursat problem
\begin{align}
\Box u^a  - q u^a &= 0, \qquad t \geq |x-a|,
\label{eq:ude}
\\
u^a(x,|x-a|) &= \frac{1}{8 \pi} \int_0^1 q(a+s(x-a)) \, ds.
\label{eq:uic}
\end{align}

We introduce some notation to state our result. For every $x \in \R^3$ and $i,j=1,2,3$, $i \neq j$, we define the angular 
derivative $\Omega_{ij} = x_i \pa_j - x_j \pa_i$.

Our main result is a uniqueness theorem for the point source inverse backscattering problem for $q$ in a special class of functions. 
Ideally, one would like to assert that if the point source backscattering data for two $q$'s are identical then the two $q$'s are identical. 
We prove this only for a special class of $q$'s. We show that two $q$'s with the same backscattering data are identical if the 
difference of the $q$'s has controlled angular derivatives as defined by (\ref{eq:angular}). 
\begin{theorem}[Uniqueness]\label{thm:uniqueness}
Suppose $q_k$, $k=1,2,$ are smooth functions on $\R^3$ with support in $B$ and $u_k^a$ are the solutions of (\ref{eq:ude}), 
(\ref{eq:uic}) for $q=q_k$ and $a \in S$. If $u_1^a(a,t) = u_2^a(a,t)$ for all $(a,t) \in S \times [0,2]$ then $q_1=q_2$ provided
there is a constant $C$, independent of $\rho$, $i,j$, such that
\beqn
\int_{|y|=\rho} | \Omega_{ij}(q_1 - q_2)(y)|^2 \, dS_y \leq C \int_{|y|=\rho} | (q_1 - q_2)(y)|^2 \, d S_y,
~ \forall \rho \in (0,1], ~ \forall ~ i,j=1,2,3.
\label{eq:angular}
\eeqn
\end{theorem}

We examine the meaning of the condition (\ref{eq:angular}).  To every $x \in \R^3$, $x \neq 0$, we associate a unique unit vector 
$\omega = x/|x|$ in $\R^3$ and a unique $\rho=|x| > 0$.
Let $\{ \phi_n(\omega) \}_{n \geq 1}$ denote an orthonormal basis for $L^2(S)$ consisting of spherical harmonics. Each 
$\phi_n(\omega)$ is the restriction to $S$ of a homogeneous harmonic polynomials $\phi_n(x)$, and the $\phi_n$ are indexed so that if 
$m<n$ then $d_m \leq d_n$ where $d_k = \text{deg} (\phi_k)$ - see [SW71]. Now $p(x)=(q_1-q_2)(x)$ has a spherical harmonic expansion
$ p(x) = \sum_{n=1}^\infty p_n(\rho) \phi_n(\omega)$ and in [RU14] it was shown that $p(x)$ satisfies the angular 
derivative condition  (\ref{eq:angular})  iff we can find $C$ (independent of $\rho$) so that
\beqn
\sum_{n=1}^\infty d_n ( d_n + 1) \, p_n(\rho)^2 
\leq C \sum_{n=1}^\infty p_n(\rho)^2, \qquad \forall \rho \in (0,1].
\label{eq:angcond}
\eeqn
Clearly (\ref{eq:angcond}) holds if $p_n(\cdot) =0$ for all $n \geq N$ for some $N$, but (\ref{eq:angcond}) also holds for some $p$ 
with infinite spherical harmonic expansions - see [RU14].

So Theorem \ref{thm:uniqueness} implies uniqueness for the inverse problem in the special case if the $q$ are radial - this was already shown in [Rak08].
- but also if the difference of the $q$'s is a finite linear combination of the spherical harmonics with 
coefficients which may be radial functions. Proving uniqueness for general $q$ is a long standing open question. There are additional 
specialized uniqueness results; if $q_1 \geq q_2$ then a uniqueness was proved in [Rak08] and if $q_1,q_2$ are small than one can 
obtain a uniqueness result analogous to Theorem 3b in [RU14] using ideas similar to those in [RU14].

The proof of Theorem \ref{thm:uniqueness}  relies on two ideas. First we need an identity obtained by using the solution of an adjoint 
problem, an idea used earlier by Santosa and Symes in [SnSy88], and by Stefanov in [St90]. The second idea used is to bound
the value of a function $f$ which is supported in $B$ by the derivatives of the mean values of $f$ on spheres centered on the boundary 
of $B$ and do it in a local manner. This bound is obtained by an idea motivated by the material in pages 185-190 in [LRS86].
%


\section{Proof of Theorem \ref{thm:uniqueness}}

Below $\cleq$ means ``less than or equal to a constant multiple of''.

\noindent
\underline{\bf Step 1}

\begin{figure}[h]
\begin{center}
\epsfig{figure=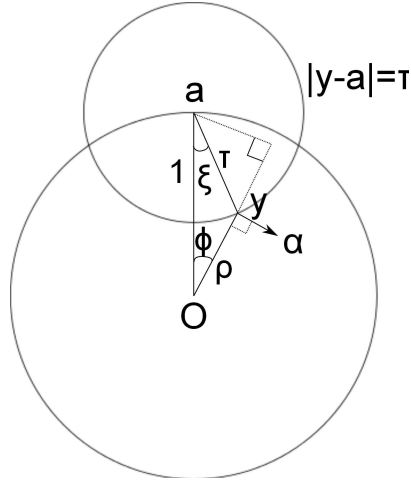, height=2.5 in}
\end{center}
\caption{Parameterization of the sphere $|y-a|=\tau$}
\label{fig:decomposition}
\end{figure}

First we estimate a function by its spherical mean values.
For any function $p(x)$ on $\R^3$ and any $\tau \geq 0$  define $(Mp)(x,\tau)$ to be the mean value of $p$
on the sphere centered at $x$ and of radius $\tau$, that is
\[
(Mp)(x,\tau) = \frac{1}{4 \pi} \int_S p(x+ \tau \omega) \, d \omega = \frac{1}{4 \pi \tau^2}
\int_{|y-x|=\tau} f(y) \, dS_y.
\]
We have the following proposition whose proof is given in subsection \ref{subsec:derivative}. 
\begin{Proposition}\label{prop:derivative}
If $p(x)$ is a smooth function on $\R^3$ with support in $B$ then, for all $(a,\tau) \in S \times (0,1)$, we have
\[
\pa_\tau \left ( \tau (Mp)(a,\tau) \right )
 = \frac{1-\tau}{2} p( (1-\tau)a ) + \frac{1}{4 \pi}
\int_{|y-a|=\tau} \frac{ (\alpha \cdot \nabla) p (y)}{\sin \phi} \, dS_y.
\]
where $\alpha(a,y)$ is a unit vector orthogonal to $y$ and $\phi$ is the angle between $y$ and $a$; see Figure \ref{fig:decomposition}.
\end{Proposition}

Let $e_i$, $i=1,2,3$ denote the standard unit vectors in $\R^3$. For any $x \in \R^3$ and for $i,j=1,2,3$ define the vectors
$T_{ij} = x_i e_j - x_j e_i$  and note that $\Omega_{ij} = T_{ij} \cdot \nabla$. For any $x,v \in \R^3$, 
we have (see Proposition 2 in [RU14])
\[
|x|^2 v = \sum_{i<j} ( v \cdot T_{ij}) T_{ij} + (v \cdot x) x.
\]
Applying this to $y, \alpha$ and noting that $\alpha$ is orthogonal to $y$ we obtain
\begin{align*}
|y|^2 ( \alpha \cdot \nabla p)(y) &= \sum_{i<j} (\alpha \cdot T_{ij}) ( T_{ij} \cdot \nabla p)(y)
=  \sum_{i<j} (\alpha \cdot T_{ij}) (\Omega_{ij} p)(y).
\end{align*}
Now $ |\alpha| =1$ and $|T_{ij}| \leq 2 |y|$, so $ |\alpha \cdot T_{ij}| \leq 2 |y|$ and hence
\[
| (\alpha \cdot \nabla p)(y)|
\leq \frac{2}{|y|} \sum_{i<j}  | (\Omega_{ij} p)(y) | .
\]
Further, the sphere $|y-a|=\tau$ may be parameterized with $\rho, \theta$ where $\rho=|y|$ and $\theta$ the rotation of $y$
about the line through $a$ (see Figure \ref{fig:decomposition}) and with this parameterization we
have\footnote
{ 
We have $\xi, \theta$
parametrizing the sphere $|y-a|=\tau$ and $dS_y = \tau^2 \sin \xi \, d \xi \, d
\theta$. Hence $dS_y = \tau^2 \sin \xi \, \dfrac{d\xi}{d \rho}  \, d \rho \, d
\theta$. The relationship between $\xi$ and $\rho$ is
$
\cos \xi = \dfrac{ 1 + \tau^2 - \rho^2}{2\tau}
$
so
$
- \sin \xi \, \dfrac{d\xi}{d\rho} = - \dfrac{\rho}{\tau},
$
hence 
$ dS_y=\tau \, \rho \, d \rho \, d \theta$.
}
$dS_y=\tau \, \rho \, d \rho \, d \theta$.
Hence Proposition \ref{prop:derivative} implies that
\begin{align}
|p((1-\tau)a)| & 
\leq \frac{2}{1-\tau} | \partial_\tau ( \tau  (Mp)(a,\tau) ) |
+ \frac{\tau}{\pi (1 - \tau)}  \sum_{i<j}  
\int_{1-\tau}^{1} 
 \int_{0}^{2 \pi}  \frac{ |\Omega_{ij} p(y)|  } {|\sin \phi|} \,  d \theta \, d \rho.
\label{eq:pP}
\end{align}

Using the law of cosines we may determine $\cos \phi$ (see Figure \ref{fig:decomposition}) and hence
\begin{align*}
\frac{1}{\sin \phi} &= \frac{2 \rho}{ \sqrt{ 4 \rho^2 - (\rho^2 + 1 - \tau^2)^2} }
=
\frac{2 \rho}{
\sqrt{ (2 \rho - (\rho^2 + 1 - \tau^2)) ( 2 \rho + (\rho^2 + 1 - \tau^2)) }
}
\\
& = \frac{2 \rho}{
\sqrt { (\tau^2 - (\rho-1)^2) ( (\rho+1)^2 -\tau^2) } } 
\\
& = \frac{2 \rho}{
\sqrt{ (\rho - (1- \tau)) ( \tau + 1 - \rho) ( \rho + 1 + \tau) ( \rho + 1 - \tau) } }.
\end{align*}
Since we will take $0< \tau \leq 1$ and $1- \tau \leq \rho \leq 1$, we have 
$\tau+1-\rho \geq \tau$, 
$\rho+1 + \tau \geq 1$,
$\rho+1-\tau \geq 1 - \tau$, 
 and we obtain
\[
\frac{1}{| \sin \phi|} \leq \frac{2\rho}{ \sqrt{\tau} \sqrt{1-\tau} \sqrt{ \rho - (1- \tau)} }.
\]
Using this in (\ref{eq:pP}) we obtain
\begin{align*}
& \quad |p( (1-\tau)a)| 
\\
&\leq  \frac{2}{1-\tau} | \partial_\tau ( \tau  (Mp)(a,\tau) ) |
+ \frac{2\sqrt{\tau} }{\pi  (1-\tau)^{3/2} } \sum_{i<j} \int_{1-\tau}^{1} 
\int_0^{2 \pi} \frac{ |\Omega_{ij} p(y)|} {\sqrt{\rho - (1-\tau)}} \, \rho \,d \theta \, d \rho.
\end{align*}
Hence, by the Cauchy-Schwarz inequality, for all $(a,\tau) \in S \times (0,1)$, we have
\begin{align}
(1-\tau)^3 & |p((1-\tau)a)|^2 
\nn
\\
& \cleq  | \partial_\tau ( \tau  (Mp)(a,\tau) ) |^2 \nn
\\
& \qquad+  \int_{1-\tau}^{1} \int_0^{2 \pi} \frac{ \rho} {\sqrt{\rho - (1-\tau)}} \, 
d \theta \, d \rho 
\,
\sum_{i<j}  \int_{1-\tau}^{1} \int_0^{2 \pi} \frac{ |\Omega_{ij} p(y)|^2} {\sqrt{\rho - (1-\tau)}} \, 
\rho \, d \theta \, d \rho
\nn
\\
& \cleq  | \partial_\tau ( \tau  (Mp)(a,\tau) ) |^2
+ \sum_{i<j}  \int_{1-\tau}^{1} \int_0^{2 \pi} \frac{ |\Omega_{ij} p(y)|^2} {\sqrt{\rho - (1-\tau)}} \, 
\rho \, d \theta \, d \rho
\nn
\\
&  \cleq  | \partial_\tau ( \tau  (Mp)(a,\tau) ) |^2
+ \sum_{i<j} \int_{|y-a|=\tau} \frac{ |\Omega_{ij} p(y)|^2} {\sqrt{|y| - (1-\tau)}} \, dS_y
\label{eq:DMp}
\end{align}
with the constant independent of $a$ and $\tau$. 
\\

\noindent
\underline{\bf Step 2}

Now suppose $q_i$, $i=1,2$ are smooth functions on $\R^3$ which are supported on $B$ and $U_i^a$ is the solution of 
(\ref{eq:Ude}), 
(\ref{eq:Uic}) when $q=q_i$. Define $v^a= U_1^a - U_2^a = u_1^a - u_2^a$ and
$p = q_1 - q_2$; then $v^a$ is the solution of the IVP
\begin{align}
\Box v^a - q_1 v^a = p U_2^a, & \qquad (x,t) \in \R^3 \times \R,
\label{eq:vde}
\\
v^a(x,t)=0, & \qquad \text{for} ~ t<0.
\label{eq:vic}
\end{align}
We also note that $v^a(x,t)$ is supported in the region $t \geq |x-a|$ and, in this region, is the smooth solution of the 
 Goursat problem
\begin{align}
\Box v^a  - q_1 v^a &= p u_2^a, \qquad t \geq |x-a|,
\label{eq:vgde}
\\
v^a(x,|x-a|) &= \frac{1}{8 \pi} \int_0^1 p(a+s(x-a)) \, ds.
\label{eq:vgic}
\end{align}
We have the following identity whose proof is given in subsection \ref{subsec:identity}.
\begin{Proposition}\label{prop:identity}
For any $\tau>0$ and all $a \in S$ we have
\begin{align}
v^a(a,2\tau) &= 
\frac{1}{8 \pi} (Mp)(a, \tau)  
+  \int_{|x-a| \leq \tau} p(x) \, k(x,\tau,a) \, dx 
\label{eq:iden}
\end{align}
where
\beqn
k(x,\tau,a) = \frac{(u_1^a+u_2^a)(x, 2\tau - |x-a|)}{4 \pi |x-a|} + \int_{|x-a|}^{2\tau-|x-a|} u_1^a(x, 2\tau-t)
\, u_2^a(x,t) \, dt, 
\label{eq:kdef}
\eeqn
\hfill if $ ~ |x-a| \leq \tau$.
\end{Proposition}
\begin{Remark} We will use a smooth extension of $k(x,t,a)$ beyond the region $|x-a| \leq t$. Note that it may seem that $k(x,t,a)$ is 
singular where $|x-a|=0$, that is at $x=a$, but $p(x)$ is zero in a neighborhood of any $a \in S$, so in (\ref{eq:iden}) we could 
replace the $k(x,t,a)$ by one with $k(x,t,a)$ multiplied by a cutoff function so that $k=0$ in a neighborhood of $x=a$.
\end{Remark}

Now, by hypothesis $v^a(a, 2 \tau) =0$ for all $(a,\tau) \in S \times [0,1]$ so,  from Proposition \ref{prop:identity}, we obtain
\[
\tau (Mp)(a, \tau) = - 8 \pi \tau \int_{|x-a| \leq \tau} p(x) \, k(x, \tau,a) \, dx, \qquad \forall a \in S, ~ \tau \in (0,1],
\]
and hence
\begin{align*}
- \pa_{\tau} ( \tau (Mp)(a,\tau) )
 = 8 \pi \tau \int_{|x-a| = \tau} & p(x) k(x, \tau,a) \, dS_x \\
& +  8 \pi
\int_{|x-a| \leq \tau} p(x) \, \pa_\tau( \tau k(x,\tau,a) ) \, dx.
\end{align*}
So
\begin{align*}
| \pa_{\tau} ( \tau (Mp)(a,\tau) ) |^2 & \cleq \int_{|x-a|=\tau} |p(x)|^2 \, dS_x
+\int_{|x-a|\leq \tau} |p(x)|^2 \, dx
\end{align*}
and hence (\ref{eq:DMp}) implies that for all $a \in S$ and all $\tau \in [0,1]$ we have
\begin{align*}
(1-\tau)^3  \, | p( (1- \tau) a) | ^2
\cleq
\int_{|y-a|= \tau} & |p(y)|^2 \, dS_y +\int_{|y-a| \leq \tau} |p(y)|^2 \, dy \\
& + \sum_{i<j}  \int_{|y-a|\leq \tau}  \frac{ |\Omega_{ij} p(y)|^2} {\sqrt{|y| - (1-\tau)}} \, dS_y.
\end{align*}
So integrating w.r.t $a$ over $S$, we obtain
\begin{align}
(1-\tau)^3  \int_S & |p((1-\tau)a)|^2 \,  dS_a 
\cleq 
\int_S \int_{|y-a|= \tau} |p(y)|^2 \, dS_y  \, dS_a
\nn \\
& + \int_S \int_{|y-a| \leq \tau} |p(y)|^2 \, dy \, dS_a
+ \sum_{i<j} \int_S \int_{|y-a|\leq \tau}  \frac{ |\Omega_{ij} p(y)|^2} {\sqrt{|y| - (1-\tau)}} \, dy
\, dS_a.
\label{eq:pop1}
\end{align}

We now convert the $dS_y \, dS_a$ integral and the $dy \, dS_a$ integral into just a $dy$ integral.
Let $e=(0,0,1)$ and let $\phi,\theta$ represent the spherical coordinates for $a$ on the sphere $S$, that is $\cos \phi = a \cdot e$. 
For any non-zero $y \in B$ and 
$\tau \in (0,1)$ we have
\[
| y - a|^2 = |y|^2 + 1 - 2 |y| \cos \phi
\]
so
\begin{align*}
\int_S \delta & ( |y-a|^2 -\tau^2) \, dS_a 
= 2 \pi \int_0^\pi \delta( |y|^2 + 1 - 2 |y| \cos \phi - \tau^2 ) \, \sin \phi \, d \phi
\\
& = 2 \pi \int_{-1}^1 \delta( |y|^2 + 1 - 2s |y| - \tau^2) \, ds
  = \frac{\pi}{|y|} \int_{-1}^1 \delta (s -  (|y|^2 + 1 - \tau^2)/(2 |y|)  ) \, ds
\\
& = \frac{\pi}{|y|} \, H( \tau + |y| -1).
\end{align*}
Also for $y \in B$, $\tau \in (0,1)$, we have
\begin{align*}
\int_S & H ( \tau^2  - |y-a|^2) \, dS_a = 2 \pi \int_0^\pi H( \tau^2 - |y|^2 - 1 + 2 |y| \cos \phi) \, \sin \phi \, d \phi 
\\
& = 2 \pi  \int_{-1}^1 H( \tau^2 - |y|^2 - 1 + 2 |y| s ) \, ds
= 2 \pi \int_{-1}^1 H ( s - (1+|y|^2 - \tau^2)/(2 |y|) ) \, ds
\\
& \leq \begin{cases} 
4 \pi, & \text{if} ~ (1+|y|^2 - \tau^2)/(2 |y|) < 1 \\
0, & \text{if} ~ (1+|y|^2 - \tau^2)/(2 |y|)  \geq 1 
\end{cases}
\\
&=4 \pi H(\tau+ |y|-1).
\end{align*}
Hence for any non-negative integrable function $f(y)$ which is supported in $B$ and any $\tau \in (0,1)$ we have
\begin{align*}
\int_S \int_{|y-a|=\tau} f(y) \, dS_y \, da
& = 2 \tau \int_S \int_{B} f(y) \, \delta(|y-a|^2-\tau^2) \, dy \, dS_a
\\
& = 2 \tau \int_{B} f(y) \int_S \delta( |y-a|^2-\tau^2) \, dS_a \, dy
\\
& =2 \pi \tau \int_{|y| \geq 1- \tau} \frac{f(y)}{|y|} \, dy 
\end{align*}
and
\begin{align*}
\int_S \int_{|y-a|\leq \tau} f(y) \, dy \, da
& = \int_S \int_{B} f(y) \, H( \tau^2 - |y-a|^2) \, dy \, dS_a
\\
& =  \int_{B} f(y) \int_S H(\tau^2 - |y-a|^2) \, dS_a \, dy
\\
& \leq 4 \pi  \int_{|y| \geq 1- \tau} f(y) \, dy.
\end{align*}

Applying these estimates to (\ref{eq:pop1}) we obtain
\begin{align}
(1-\tau)^3  \int_S |p((1-\tau)a)|^2 \, dS_a 
& \cleq 
\int_{|y|\geq 1- \tau} \frac{|p(y)|^2}{|y|} dy
+ \sum_{i<j}  \int_{|y|\geq 1-\tau}  \frac{ |\Omega_{ij} p(y)|^2} {\sqrt{|y| - (1-\tau)}} dy.
\label{eq:pop2}
\end{align}
Define
\[
P(\rho) := \int_{|y|=\rho} |p(y)|^2 \, dS_y;
\]
then the angular control condition (\ref{eq:angular}) is equivalent to the statement that
\[
\sum_{i<j}  \int_{|y|=\rho}  |\Omega_{ij} p(y)|^2\, dS_y
\cleq P(\rho), \qquad \forall \rho \in (0,1).
\]
Hence (\ref{eq:pop2}) implies that for all $\tau \in (0,1)$ we have
\[
(1-\tau) P(1-\tau) \cleq \int_{1-\tau}^1 \frac{P(\rho)}{\rho} \, d\rho + \int_{1-\tau}^1 \frac{ P(\rho)}{\sqrt{\rho - (1-\tau)} } \, d \rho.
\]
which may be rewritten as
\[
s P(s) \cleq \int_{s}^1 \frac{P(\rho)}{\rho} \, d\rho + \int_{s}^1 \frac{ P(\rho)}{\sqrt{\rho - s} } \, d \rho
\]
for all $s \in (0,1)$. We fix an $\ep>0$, then for all $s \in (\ep,1)$ we have a constant $C_\ep$ so that
\[
P(s) \leq C_\ep \int_s^1  \frac{ P(\rho)}{ \sqrt{\rho -s}} \, d \rho, \qquad \forall s \in [\ep,1].
\]
Reusing this inequality in the relation we obtain, for all $s \in [\ep,1]$,
\begin{align*}
P(s) &\leq C_\ep^2 \int_s^1 \int_\rho^1 \frac{ P(r)}{ \sqrt{ (\rho-s) (r-\rho)}} \, dr \, d \rho
= C_\ep^2 \int_s^1 P(r) \int_s^r \frac{ 1}{ \sqrt{ (\rho-s) (r-\rho)}} \, d \rho \, d r
\\
&= \pi C_\ep^2 \int_s^1 P(r) \, dr.
\end{align*}
Hence by Gronwall's inequality, $P(s)=0$ for all $s \in [\ep,1]$.


\section{Proofs of propositions}

\subsection{Proof of Proposition \ref{prop:derivative}}\label{subsec:derivative}

We have 
\begin{align*}
(Mp)(a,\tau) & = \frac{1}{4 \pi \tau^2} \int_{|y-a|=\tau} p(y) \, dS_y
= \frac{1}{4 \pi \tau^3} \int_{|y-a| \leq \tau} \nabla \cdot ( (y-a) \, p(y) ) \, dy,
\end{align*}
hence
\begin{align*}
\frac{d}{d \tau} ( 4 \pi \tau^3 (Mp)(a, \tau) )
&= \int_{|y-a|=\tau}  \nabla \cdot ( (y-a) \, p(y) ) \, dS_y,
\\
&= 3 \int_{|y-a|}  p(y) \, dS_y + \int_{|y-a|=\tau} (y-a) \cdot \nabla p(y) \, dS_y
\\
& = 12 \pi \tau^2 (Mp)(a, \tau) + \int_{|y-a|=\tau} (y-a) \cdot \nabla p(y) \, dS_y
\end{align*}
implying
\[
4 \pi \tau^3 \frac{d}{d \tau} (Mp)(a, \tau) =\int_{|y-a|=\tau} (y-a) \cdot \nabla p(y) \, dS_y
\]
so
\begin{align}
\pa_\tau( \tau (Mp)(a,\tau)) &= \tau \pa_\tau( (Mp)(a,\tau)) + (Mp)(a,\tau)
\nn
\\
& = (Mp)(a,\tau) + \frac{1}{4 \pi \tau^2} \int_{|y-a|=\tau} (y-a) \cdot \nabla p(y) \, dS_y  .
\label{eq:dMp}
\end{align}

To any $y \in \R^3$, we associate its spherical coordinates $(\rho, \theta, \phi)$ where $\rho =|y|$,
$\phi$ is the angle between $y$ and $a$, and $\theta$ is the rotation about $a$. Also, let $\alpha$ be the downward pointing unit vector, in the plane containing $y$ and $a$, so that $\alpha \perp y$. See Figure \ref{fig:decomposition}.

First we decompose $y-a$ as the sum of a vector parallel to $y$ and a vector parallel to $\alpha$.
From the right triangle we see that the length of the projection of $y-a$ on $y$ is $\cos \phi - \rho$ and the length of the component of $y-a$ orthogonal to $y$ is $\sin \phi$. Hence
\[
y-a = (\rho - \cos \phi) \, \frac{y}{\rho} + \sin \phi \, \alpha
\]
so
\beqn
(y-a) \cdot (\nabla p)(y) = (\rho - \cos \phi) \frac{\pa p}{\pa \rho}
+ \sin \phi \, (\alpha \cdot \nabla) p (y). 
\label{eq:ttemp1}
\eeqn

The $ \dfrac{\pa p}{\pa \rho}$ in (\ref{eq:ttemp1}) is the partial derivative of $p(\rho,\phi, \theta)$ when $\rho, \phi, \theta$ are regarded as independent variables. However, we have parameterized the sphere $|y-a|=\tau$ with $\rho,\theta$ (and $\tau$) and 
$\phi$ is regarded as function of $\rho$ via the relation
\beqn
\cos \phi = \frac{\rho^2 + 1 - \tau^2}{2 \rho}
\label{eq:phirho}
\eeqn
which implies 
\[
2 \cos \phi - 2 \rho \sin \phi \, \frac{ d \phi}{d \rho} = 2 \rho.
\]
With this we observe that
\begin{align*}
\frac{\pa}{\pa \rho} ( p( \rho, \phi(\rho), \theta) ) = \frac{\pa p}{\pa \rho} + 
\frac{\pa p}{\pa \phi} \, \frac{\pa \phi}{\pa \rho}
=  \frac{\pa p}{\pa \rho} + \frac{\pa p}{\pa \phi} \, \frac{\cos \phi - \rho}{\rho \sin \phi}.
\end{align*}
Now, at $y$, as $\phi$ increases by $d \phi$, $y$ moves along a circular arc of radius $|y|$ through a distance $|y| \, d \phi$, hence
\[
\frac{\pa p }{ \pa \phi} = |y| \, (\alpha \cdot \nabla p)(y) = \rho  \, (\alpha \cdot \nabla p)(y)
\]
which implies that
\[
\frac{\pa p}{\pa \rho} = \frac{\pa}{\pa \rho} ( p( \rho, \phi(\rho), \theta) ) + \frac{ \rho - \cos \phi}{ \sin \phi}
 (\alpha \cdot \nabla p)(y).
\]
Hence from (\ref{eq:ttemp1}) and (\ref{eq:phirho})
\begin{align*}
& (y-a) \cdot (\nabla p)(y) 
\\
& =  (\rho - \cos \phi) \, \frac{\pa}{\pa \rho} ( p( \rho, \phi(\rho), \theta) )
+  \frac{ (\rho - \cos \phi)^2}{ \sin \phi}
 (\alpha \cdot \nabla p)(y)
 +  \sin \phi \, (\alpha \cdot \nabla) p (y)
 \\
 &= (\rho - \cos \phi) \, \frac{\pa}{\pa \rho} ( p( \rho, \phi(\rho), \theta) ) + 
 \frac{ \rho^2 + 1 - 2 \rho \cos \phi}{\sin \phi} \,  (\alpha \cdot \nabla) p (y)
 \\
 & =  \frac{\rho^2 + \tau^2 -1}{2 \rho} \, \frac{\pa}{\pa \rho} ( p( \rho, \phi(\rho), \theta) ) + 
 \frac{\tau^2}{\sin \phi}  (\alpha \cdot \nabla) p (y)
\end{align*}

With the $\rho, \theta$ parameterization of the sphere $|y-a|=\tau$ we 
have\footnote
{ 
We have $\xi, \theta$
parametrizing the sphere $|y-a|=\tau$ and $dS_y = \tau^2 \sin \xi \, d \xi \, d
\theta$. Hence $dS_y = \tau^2 \sin \xi \, \dfrac{d\xi}{d \rho}  \, d \rho \, d
\theta$. The relationship between $\xi$ and $\rho$ is
$
\cos \xi = \dfrac{ 1 + \tau^2 - \rho^2}{2\tau}
$
so
$
- \sin \xi \, \dfrac{d\xi}{d\rho} = - \dfrac{\rho}{\tau},
$
hence 
$ dS_y=\tau \, \rho \, d \rho \, d \theta$.
}
$dS_y=\tau \, \rho \, d \rho \, d \theta$ so
\begin{align*}
\int_{|y-a|=\tau} (y-a) \cdot (\nabla p)(y)
& = \frac{\tau}{2} \int_0^{2 \pi} \int_{1-\tau}^{1+\tau}  (\rho^2 + \tau^2 -1) \, \frac{\pa }{\pa \rho} ( p(\rho, \phi(\rho),\theta) )
\, d \rho \, d \theta
\\
& \qquad 
+ \tau^2 \int_{|y-a|=\tau} \frac{ (\alpha \cdot \nabla) p (y)}{\sin \phi} \, dS_y
\\
&=2 \pi \tau^2 (1- \tau) ( p((1-\tau)a)
- \tau\int_0^{2 \pi} \int_{1-\tau}^{1+\tau} p(y) \, \rho \, d \rho \, d \theta
\\
& \qquad
+ \tau^2 \int_{|y-a|=\tau} \frac{ (\alpha \cdot \nabla) p (y)}{\sin \phi} \, dS_y
\\
& = 2 \pi \tau^2 (1- \tau)  p((1-\tau)a) - \int_{|y-a|=\tau} p(y) \, dS_y
\\
& \qquad
+ \tau^2 \int_{|y-a|=\tau} \frac{ (\alpha \cdot \nabla) p (y)}{\sin \phi} \, dS_y.
\end{align*}
Hence from (\ref{eq:dMp}) we obtain
\begin{align}
\pa_\tau \left ( \tau (Mp)(a,\tau) \right )
& = \frac{1-\tau}{2} p( (1-\tau)a ) + \frac{1}{4 \pi}
\int_{|y-a|=\tau} \frac{ (\alpha \cdot \nabla) p (y)}{\sin \phi} \, dS_y.
\label{eq:rigor}
\end{align}

\subsection{Proof of Proposition \ref{prop:identity}}\label{subsec:identity}

We give a formal proof of the proposition, based on (\ref{eq:Ude}), (\ref{eq:Uic}) and (\ref{eq:vde}),
(\ref{eq:vic}). A rigorous proof can be constructed using (\ref{eq:ude}), (\ref{eq:uic}) and
(\ref{eq:vgde}), (\ref{eq:vgic}). We also note that if $M$ is the surface $\phi(z)=0$ in $\R^n$ then
\beqn
\int_M h(z) \, dS_z = \int_{\R^n} h(z) \, |\nabla \phi(z)| \, \delta( \phi(z) ) \, dz.
\label{eq:deltaS}
\eeqn

We note that $U_1^a(x, 2\tau-t)$ is the solution of the final value problem
\begin{align*}
(\Box - q_1) U_1^a(x, 2\tau -t) = \delta(x-a, 2\tau-t), & \qquad (x,t) \in \R^3 \times \R
\\
U_1^a(x, 2\tau -t) = 0, & \qquad t > 2\tau.
\end{align*}
Further, $v(x,t)=0$ for $t<0$, and for each fixed $t$, $U_2(x,t)$ has compact support in $x$. So using 
(\ref{eq:vde}) and an integration by parts we have
\begin{align}
\int_{\R^3 \times \R}   p(x) \, U_2^a(x,t) \, U_1^a(2\tau -t, x) \, dx \, dt
&= \int_{\R^3 \times \R} (\Box - q_1) v^a(x,t) \, U_1^a(2\tau -t, x) \, dx \, dt
\nn
\\
& =  \int_{\R^3 \times \R} v^a(x,t) \, (\Box - q_1) U_1^a(2\tau -t, x) \, dx \, dt
\nn
\\
& = \int_{\R^3 \times \R} v^a(x,t) \, \delta( x-a, 2\tau-t) \, dx \, dt
\nn
\\
& = v^a(a,2\tau).
\label{eq:UUv}
\end{align}
On the other hand
\begin{align*}
& \quad \int_{\R^3 \times \R}  p(x) \, U_2^a(x,t) \, U_1^a(2\tau -t, x) \, dx \, dt
\\
 & = \int_{\R^3 \times \R} p(x) \left ( \frac{ \delta(t - |x-a|)}{4 \pi |x-a|} + u_2^a(x,t) \right )
\left ( \frac{ \delta(2\tau -t - |x-a|)}{4 \pi |x-a|} + u_1^a(x,2\tau-t) \right ) dx dt
\\
&=\frac{1}{16 \pi^2} \int_{\R^3 \times \R} \frac{ p(x) \, \delta(t - |x-a|) \, \delta(2\tau - t - |x-a|)}{|x-a|^2}
 \, dx \, dt
\\
& \qquad +
\frac{1}{4 \pi} \int_{\R^3 \times \R} \frac{p(x)}{|x-a|} \, \delta(t - |x-a|) \, u_1^a(x, 2\tau-t) \, dx \, dt
\\
& \qquad +
\frac{1}{4 \pi} \int_{\R^3 \times \R} \frac{ p(x)}{|x-a|} \,  \delta(2\tau-t - |x-a|) \, u_2^a(x, t) \, dx \, dt
\\
& \qquad +
 \int_{\R^3 \times \R} p(x) \, u_1^a(x, \tau-t) \, u_2^a(x,t) \, dx \, dt
\\
&= \frac{1}{32 \pi^2 \tau^2} \int_{|x-a|=\tau} p(x) \, dS_x 
+ \frac{1}{4 \pi} \int_{\R^3} p(x) \, \frac{ ( u_1^a + u_2^a)(x, 2\tau - |x-a|) }{|x-a|}\, dx 
\\
&  \qquad + \int_{\R^3 \times \R} p(x) \, u_1^a(x,2 \tau-t) \, u_2^a(x,t) \, dx \, dt.
\end{align*}
Since $u_i^a(x, 2\tau-|x-a|)$ is supported in $ |x-a| \leq 2\tau - |x-a|$, that is $|x-a| \leq \tau$, the second integral is over the region $|x-a| \leq \tau$. Since $u_1^a(x, 2\tau-t)$ is supported in the downward conical region $|x-a| \leq 2\tau -t$ and $u_2^a(x,t)$ is supported in the upper conical region $|x-a| \leq t$, the product of these two will be supported in the double conical region
\[
\{ (x,t) :  |x-a| \leq t \leq 2\tau - |x-a|, ~ |x-a| \leq \tau \}.
\]
Hence
\begin{align*}
\int_{\R^3 \times \R} &  p(x) \, U_2^a(x,t) \, U_1^a(2\tau -t, x) \, dx \, dt
\\
& = \frac{1}{32 \pi^2 \tau^2} \int_{|x-a|=\tau} p(x) \, dS_x 
+ \frac{1}{4 \pi} \int_{|x-a| \leq \tau} p(x) \, \frac{( u_1^a + u_2^a)(x, 2\tau - |x-a|)}{|x-a|} \, dx 
\\
&  \qquad + \int_{|x-a| \leq \tau} p(x)  \left ( \int_{|x-a|}^{ 2\tau - |x-a|}  u_1^a(x,2 \tau-t) \, u_2^a(x,t) \, dt 
\right ) \, dx.
\\
\end{align*}
So using (\ref{eq:UUv}) we obtain Proposition \ref{prop:identity}.
%


\bibliographystyle{amsalpha}

\end{document}